\documentclass[a4paper, 12pt]{article}
\usepackage{amsmath}
\usepackage{amsthm}
\usepackage[all]{xy}
\usepackage{amsfonts}

\title{Barr's Embedding Theorem for Enriched Categories}
\author{Dimitri Chikhladze}

\begin{document}
\newtheorem{definition}{Definition}
\newtheorem{proposition}[definition]{Proposition}
\newtheorem{lemma}[definition]{Lemma}
\newtheorem{theorem}[definition]{Theorem}
\maketitle
\abstract{We generalize Barr's embedding theorem for regular categories to the context of enriched categories.}

\section{Introduction}
M. Barr proved that a (small) regular category can be embedded into a category of (small) presheaves in \cite{Barr1}. Barr gave a different proof of the same result in \cite{Barr}, and other proofs have appeared (F. Borceux, \cite{Borceux}, M. Makkai, \cite{Makkai}). Our purpose is to generalize the theorem from the case of ordinary (set enriched) categories to that of categories enriched in a monoidal category. We are influenced by \cite{Barr}, although the paper contains some inaccuracies, which the author acknowledged and outlined how to fix them in correspondence.

We adopt a notion of regularity for enriched categories suitable for our setting. Regular enriched categories have been considered before by B. Day, R. Street in \cite{Street}. 

\section{Regular categories}
Suppose $\mathbb{C}$ is a $\mathbb{V}$-category. By a morphism in $\mathbb{C}$ we mean a morphism in the underlying ordinary category $\mathbb{C}_0$. Speaking of (co)limits of diagrams of such morphisms we mean conical (co)limits in $\mathbb{C}$ (these, by definition, are preserved by each of the functors $\mathbb{C}(-, A)$ ($\mathbb{C}(A, -)$), a condition stronger than just being a colimit in the underlying ordinary category). Particularly we can speak of filtered colimits in a $\mathbb{V}$-category. For general theory of enriched categories we refer to G.M. Kelly \cite{KellyBook}. 

Suppose $\mathbb{V}$ is a locally finitely presentable symmetric monoidal closed category. In particular, filtered colimits in $\mathbb{V}$ commute with finite limits. 

Let $\mathbb{L}$ be a locally finitely presentable $\mathbb{V}$-category. This means that $\mathbb{L} = Lex(\mathbb{C}, \mathbb{V})$, where $\mathbb{C}$ is a finitely complete $\mathbb{V}$-category, i.e. $\mathbb{C}$ admits finite conical limits and cotensor products with finitely presentable objects of $\mathbb{V}$; objects of $\mathbb{L}$ are the functors $F: \mathbb{C} \rightarrow \mathbb{V}$ which preserve these finite limits. There is a duality between $\mathbb{V}$-categories with finite limits and locally finitely presentable $\mathbb{V}$-categories, whereby a locally finitely presentable $\mathbb{V}$-category corresponds to a $\mathbb{V}$-category with finite limits as described above, and conversely the subcategory of the finite objects of a locally finite presentable $\mathbb{V}$-category is the $\mathbb{V}$-category with finite limits corresponding to it. In the literature this duality is known as Gabriel-Ulmer duality. For more on finitely presentable categories see P. Gabriel, F. Ulmer \cite{GU} and G.M. Kelly \cite{Kelly}. The former considers only the case $\mathbb{V} = Set$, while the latter deals with the case of a general $\mathbb{V}$ still being essentially self-contained. 

In the next proposition we include some, mostly well known, properties of $\mathbb{L}$ we need.

\begin{proposition} \label{prop}
The following hold.
\begin{enumerate}
\item[(1)] $\mathbb{L}$ is complete and cocomplete.
\item[(2)] The Yoneda embedding $Y: \mathbb{C} \rightarrow \mathbb{L}^{op}$ preserves colimits and finite limits.
\item[(3)] Filtered colimits in $\mathbb{L}$ commute with finite limits, in particular a filtered colimit of regular monos (see below) is a regular mono.
\item[(4)] For every representable presheaf $A$, the functor $\mathbb{L}(A, -)$ commutes with filtered colimits.
\end{enumerate}
\end{proposition}

\begin{proof}
It is well known that the inclusion $\mathbb{L} \rightarrow [\mathbb{C}, \mathbb{V}]$ has a left adjoint. This implies the cocompletness. It also implies completeness since a limit of left exact presheaves in the presheaf category is left exact. We have (1). 

Clearly, $Y: \mathbb{C} \rightarrow \mathbb{L}^{op}$ will preserve all those limits in $\mathbb{C}$ which are preserved by each presheaf in $\mathbb{L}$. Thus, it preserves finite limits. The Yoneda embedding preserves colimits. We have (2).

In a presheaf category limits and colimits are pointwise. So, if filtered colimits commute with finite limits in $\mathbb{V}$, also in $[\mathbb{C}, \mathbb{V}]$ filtered colimits commute with finite limits. Now to prove (3) we only need to see that the inclusion of $\mathbb{L}$ into $[\mathbb{C}, \mathbb{V}]$ preserves filtered colimits and finite limits. The preservation of filtered colimits follows from the fact that in the presheaf category a filtered colimit of left exact presheaves is left exact, which is a consequence of exactness of filtered colimits in $\mathbb{V}$. Being a right adjoint the inclusion preserves limits.

Similarly, representables commute with all colimits in the presheaf category. Hence representables commute with all those colimits in $\mathbb{L}$ which are preserved by the inclusion of $\mathbb{L}$ into $[\mathbb{C}, \mathbb{V}]$. We have already observed that the inclusion preserves filtered colimits, hence (4).

\end{proof}

A morphism of $\mathbb{C}$ is a regular epi if it is a coequalizer. If a morphism is a coequalizer, then it is the coequalizer of its kernel pair if the latter exists. So then, a morphism $e: C \rightarrow D$ is a regular epi iff it is the coequalizer of its kernel pair. We will say that a regular epi is stable under pullbacks if a pullback of it by any arrow is again a regular epi. We will say that a regular epi is stable under cotensors with the finite objects if it remains a regular epi after cotensoring with any finite object of $\mathbb{V}$. We consider the following notion of regularity for enriched categories.

\begin{definition}
A $\mathbb{V}$-category $\mathbb{C}$ is regular if it is finitely complete, coequalizers of kernel pairs exist and the regular epis are stable under pullbacks and cotensors with the finite objects. A functor is regular if it preserves finite limits and regular epis.
\end{definition}

Below we give necessary and sufficient conditions for regularity of a $\mathbb{V}$-category $\mathbb{C}$ in a way to make clear the relationship with the (ordinary) regularity of the underlying ordinary category $\mathbb{C}_0$; all concepts (regularity, finite limits etc.) in these conditions relating to $\mathbb{C}_0$ are in the ordinary sence. A $\mathbb{V}$-category $\mathbb{C}$ is regular if and only if the following hold.

\begin{enumerate}
	\item The underlying ordinary category $\mathbb{C}_0$ is regular.
	\item The finite limits in $\mathbb{C}_0$ are preserved by each of the functors $\mathbb{C}(C, -) : \mathbb{C}_0 \rightarrow \mathbb{V}$.
	\item The coequalizers of kernel pairs in $\mathbb{C}_0$ are preserved by each of the functors $\mathbb{C}(-, C) : \mathbb{C}_0 \rightarrow \mathbb{V}$.
	\item $\mathbb{C}$ admits cotensors by the finite objects of $\mathbb{V}$; a cotensor product of an object $C$ of $\mathbb{C}$ with a finite object $V$ of $\mathbb{V}$ is written as $C^V$.
	\item If a map $e : C\rightarrow D$ is a regular epi in $\mathbb{C}_0$, then for any finite object $V$ of $\mathbb{V}$ the map $e^V : C^V\rightarrow D^V$ is a regular epi.
\end{enumerate}

Much like the set enriched case, a $\mathbb{V}$-enriched regular category admits regular epi-mono factorization system. In fact a definition of regularity involving this factorization can be given. It is instructive to note that a functor is regular iff it preserves finite limits and coequalizers of kernel pairs. 

The dual notion of regularity is the notion of coregularity. Under this duality, regular epis correspond to regular monos; pullback stability becomes pushout stability; stability under cotensor products with the finite objects becomes stability under tensor product with the finite objects; and regular epi-mono factorization becomes epi-regular mono factorization.

The following theorem asserts that if a $\mathbb{V}$-category with finite limits is regular then the finitely presentable $\mathbb{V}$-category corresponding to it through Gabriel-Ulmer duality is coregular.

\begin{theorem} \label{proposition}
If $\mathbb{C}$ is regular, then $Lex(\mathbb{C}, \mathbb{V})$ is coregular. 
\end{theorem}

\begin{proof}
Since $\mathbb{C}$ is regular $\mathbb{C}^{op}$ is coregular. So, $\mathbb{L}$ contains a coregular subcategory of the representable presheaves, which we can identify with $\mathbb{C}^{op}$. Throughout the proof we will routinely use the fact that the embedding $\mathbb{C}^{op} \rightarrow \mathbb{L}$ preserves finite colimits and limits.

Let us prove that the pushout of a regular mono is a regular mono. 

Lemma in \cite{Street} proves that each diagram in $\mathbb{L}$ the indexing category of which has finite homsets can be written as a filtered colimit of diagrams landing in the subcategory of representables. An instance of this lemma is that each pushout diagram in $\mathbb{L}$ can be written as a colimit of representable pushout diagrams. We will show that a pushout diagram in $\mathbb{L}$ with one specified arrow a regular mono can be written as a colimit of representable pushout diagrams in all of which the specified arrow is a regular mono. The set-based case of this fact appears in \cite{Barr}. 

Let $K$ be the graph $1 \leftarrow 2 \rightarrow 3$ so that $q: K \rightarrow \mathbb{L}_0$ is a diagram in $\mathbb{L}$ of the form

\bigskip

$\stackrel{q_1}{\longleftarrow} \stackrel{q_2}{\longrightarrow}$. 

\bigskip

\noindent We fix a pushout diagram $u: K \rightarrow \mathbb{L}_0$ in $\mathbb{L}$, and take $\mathbb{D}$ to be the comma category

\bigskip

$\xy
(0,0)*{}="A"; 
(22,0)*{}="B";
(0,-18)*{}="C";
(22, -18)*{}="D";
(11,-9)*{\Rightarrow}; 
"A"*{\mathbb{D}}; 
"B"*{1}; 
"C"*{[K, \mathbb{C}_0^{op}]};
"D"*{[K, \mathbb{L}_0]};  
{\ar^{} "A"+(3,0); "B"+(-2,0)}
{\ar_{[K, Y_0]} "C"+(6,0); "D"+(-6,0)}
{\ar_{r} "A"+(-2,-3); "C"+(-2,3)}
{\ar^{u} "B"+(-2,-3); "D"+(-2,3)}
\endxy
$

\bigskip

\noindent By the Lemma in \cite{Street}, $\mathbb{D}$ is filtered and  

\bigskip

$colim([K, Y_0]r) = u$. 

\bigskip

\noindent An object of $\mathbb{D}$ consists of a graph $q: K \rightarrow \mathbb{L}_0$ landing in the subcategory of representable presheaves and a natural transformation $t: q \rightarrow u$. 

\bigskip

$
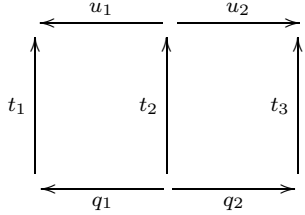
\xy
(0,0)*{}="A"; 
(16,0)*{}="B";
(32,0)*{}="C";
(0,-18)*{}="D"; 
(16,-18)*{}="E";
(32,-18)*{}="F";  
{\ar_{u_1} "B"+(2,2); "A"+(2,2)}
{\ar^{u_2} "B"+(4,2); "C"+(4,2)}
{\ar^{t_1} "D"; "A"}
{\ar^{t_2} "E"; "B"}
{\ar^{t_3} "F"; "C"}
{\ar^{q_1} "E"+(-3,-2); "D"+(-3,-2)}
{\ar_{q_2} "E"+(-3,-2); "F"+(-3,-2)}
\endxy
$

\bigskip
\bigskip

\noindent Also $t$ is the coprojection into the colimit $u$ at $q$.

Suppose now that $u_2$ is a regular mono. Take $\mathbb{D}'$ to be a subcategory of $\mathbb{D}$ consisting of those objects $q = (q, t)$ for which $q_2$ is a regular mono. We will show that the inclusion of $\mathbb{D}'$ into $\mathbb{D}$ is a final functor.  

Let $q = (q, t)$ be any object of $\mathbb{D}$. Since the subcategory of the representables is coregular every morphism in it can be uniquely factored into an epi followed by a regular mono. Let $q_2 = q'_2e$ be this unique factorization for $q_2$. Regular monos have the diagonal fill in property with respect to epis hence there exists a morphism $t'_2$ such that $t'_2e = t_2$ and $u_2t'_2 = t_3q'_2$. Let $q'_1$ be the pushout of $q_1$ by $e$. We can assume that $q'_1$ is a morphism between representables.

\bigskip

$
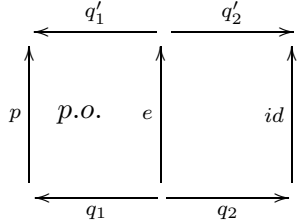
\xy
(0,0)*{}="A"; 
(16,0)*{}="B";
(32,0)*{}="C";
(0,-18)*{}="D"; 
(16,-18)*{}="E";
(32,-18)*{}="F"; 
(8,-9)*{p.o.}; 
{\ar_{q'_1} "B"+(2,2); "A"+(2,2)}
{\ar^{q'_2} "B"+(4,2); "C"+(4,2)}
{\ar^{p} "D"; "A"}
{\ar^{e} "E"; "B"}
{\ar^{id} "F"; "C"}
{\ar^{q_1} "E"+(-3,-2); "D"+(-3,-2)}
{\ar_{q_2} "E"+(-3,-2); "F"+(-3,-2)}
\endxy
$

\bigskip
\bigskip

\noindent There exists a unique $t'_1$ with $t'_1p = t_1$ and $t'_1q'_1 = u_1t'_2$. 

For any object $(q, t)$ of $\mathbb{D}$ we obtained an object $(q', t')$ of $\mathbb{D'}$, with the triple $(p, e, id)$ becoming a morphism $(q, t) \rightarrow (q', t')$ in $\mathbb{D}$. Moreover, any other morphism from $(q, t)$ to an object of $\mathbb{D'}$ factors through $(p, e, id)$. From here we can infer that the function taking $(q, t)$ to $(q', t')$ is an object function for a functor $\mathbb{D} \rightarrow \mathbb{D'}$ which is a left adjoint to the inclusion $i: \mathbb{D'} \rightarrow \mathbb{D}$. As a consequence the inclusion $i$ is a final functor. Hence we have:

\bigskip

$colim([K, Y_0]ri) = u$.

\bigskip

Now we prove that a pushout of a regular mono is a regular mono by showing this for our generic pushout diagram $u$. Let $q^\star_1 : [K, \mathbb{L}_0] \rightarrow [\rightarrow, \mathbb{L}]$ be the functor which takes $q$ to the pushout of $q_2$ by $q_1$. Because of coregularity of the subcategory of the representables, this functor sends the objects in the image of $[K, Y_0]ri: D' \rightarrow [K, \mathbb{L}_0]$ to regular monos. The following calculation expresses $q^\star_1(u)$ as a colimit of regular monos:      

\bigskip

$q^\star_1(u) = q^\star_1(colim([K, Y_0]ri)) = colim(q^\star_1([K, Y_0]ri))$.

\bigskip

\noindent By Proposition \ref{prop} (3), $q^\star_1(u)$ is a regular mono.

Stability under tensors with the finite objects follows from the proposition below, which provides a more general fact about $\mathbb{L}$.
\end{proof}

\begin{proposition} \label{prop3}
If regular epis in $\mathbb{C}$ are stable under cotensors with the finite objects, then for any object $V$ of $\mathbb{V}$ the functor $V\cdot - : \mathbb{L} \rightarrow \mathbb{L}$ preserves regular monos.
\end{proposition}

\begin{proof}
Every object $V$ of $\mathbb{V}$ is a filtered colimit of finite objects. The functor $-\cdot L : \mathbb{V} \rightarrow \mathbb{L}$ preserves colimits because it has a right adjoint the functor $\mathbb{L}(L, -) : \mathbb{L} \rightarrow \mathbb{V}$. So, given a morphism $m$ in $\mathbb{L}$, $V\cdot m$ is a filtered colimit of morphisms of the form $U\cdot m$ with $U$ a finite object. Since a filtered colimit of regular monos is a regular mono, by Proposition \ref{prop} (3), $V\cdot m$ is a regular mono if each $U\cdot m$ is a regular mono. This reduces proving the proposition to proving that regular monos are stable under tensoring with the finite objects. After the argument in the previous proof, it should be obvious that any regular mono $m$ can be written as a filtered colimit of regular monos between representables. Thus further, we only need to show the stability under tensors with the finite objects of the regular monos between representables. This is clearly true under the assumption of the proposition, using the fact again that the subcategory of the representables is equivalent to $\mathbb{C}^{op}$ and its embedding into $\mathbb{L}$ preserves tensors with finite objects and regular monos.  
\end{proof}

\section{Embedding theorem} 
Suppose that $\mathbb{C}$ is a regular $\mathbb{V}$-category.  

\begin{lemma} \label{lemma}
For every object $F$ of $\mathbb{L}$, there exists a regular mono $q: F \rightarrow F_\infty$ in $\mathbb{L}$ such that for each regular mono $a: A \rightarrow B$ between representables there exists a map $v: \mathbb{L}(A, F) \rightarrow \mathbb{L}(B, F_\infty)$ for which the diagram

\bigskip

$\xy
(0,0)*{}="A"; 
(36,0)*{}="B";
(18,-18)*{}="C";
"A"*{\mathbb{L}(A, F)};
"B"*{\mathbb{L}(A, F_\infty)};  
"C"*{\mathbb{L}(B, F_\infty)};    
{\ar^{q_\ast} "A"+(8,0); "B"+(-9,0)}
{\ar_{v} "A"+(3,-3); "C"+(-3,3)}
{\ar_{a^\ast} "C"+(3,3); "B"+(-3,-3)}
\endxy
$

\bigskip 
 
\noindent is commutative.

\end{lemma}

\begin{proof}
Choose a well order on the set of all regular monos $a: A \rightarrow B$ of representables. By transfinite induction on this well ordered set define $F_a$ as a pushout of

\bigskip

$\mathbb{L}(A, F) \cdot B \stackrel{1\cdot a}{\longleftarrow} \mathbb{L}(A, F)\cdot A \stackrel{i_a}{\longrightarrow} F'_a$

\bigskip

\noindent in $\mathbb{L}$, where $F'_0 = F$ and $i_0$ is the evaluation for the first ordinal, and $F'_a = colim_{b < a} F_b$ and $i_a$ is the evaluation followed by the obvious map $F \rightarrow F'_a$ for other ordinals. 

Take $F_\infty$ to be $colim_{a} F_a$. All the obvious maps $F \rightarrow F_a$ are regular monos since regular monos are pushout and filtered colimit stable in $\mathbb{L}$. Take $q$ to be the map $F \rightarrow F_\infty$, a regular mono too. To check the required property, observe that $F_\infty$ also is a colimit of the diagram 

\bigskip

$\mathbb{L}(A, F) \cdot B \stackrel{1\cdot a}{\longleftarrow} \mathbb{L}(A, F)\cdot A \stackrel{ev}{\longrightarrow} F$

\bigskip

\noindent in which $a: A \rightarrow B$ varies over all regular monos between all representables. For a regular $a$, the needed $v: \mathbb{L}(A, F) \rightarrow \mathbb{L}(B, F_\infty)$ is determined by transposing $B$ in the coprojection $v: \mathbb{L}(A, F)\cdot B \rightarrow F_\infty$. The commutativity of the triangle is straightforward to verify. 
\end{proof}

\begin{proposition}
For any left exact functor $F : \mathbb{C} \rightarrow \mathbb{V}$ there exists a regular functor $P : \mathbb{C} \rightarrow \mathbb{V}$ and a regular mono $F \rightarrow P$.
\end{proposition}

\begin{proof}
Define $P$ to be a colimit of a diagram

\bigskip

$F \longrightarrow F_1 \longrightarrow F_2 \longrightarrow \ldots$  

\bigskip

\noindent in $\mathbb{L}$, where $F_{n+1} = (F_{n})_\infty$ are given by the last lemma. The canonical arrow $F \rightarrow P$ is a regular mono. We should show that $P$ is a regular functor.

The following observation trivially follows from the Yoneda Lemma. A functor $P$ preserves regular epis iff for every regular mono $a: A \rightarrow B$ between representables, $a^\ast : \mathbb{L}(B, P) \rightarrow \mathbb{L}(A, P)$ is a regular epi. Back to our setting, given such a regular mono $a$ there is a diagram

\bigskip

$\xy
(0,0)*{}="A"; 
(24,0)*{}="B";
(48,0)*{}="C";
(0,-18)*{}="D"; 
(24,-18)*{}="E";
(48,-18)*{}="F";
"A"*{\mathbb{L}(A, F)}; 
"B"*{\mathbb{L}(A , F_1)}; 
"C"*{\mathbb{L}(A, F_2)}; 
"D"*{\mathbb{L}(B, F)}; 
"E"*{\mathbb{L}(B , F_1)}; 
"F"*{\mathbb{L}(B, F_2)};
(62, 0)*{\ldots};
(62, -18)*{\ldots};  
{\ar_{} "A"+(8,0); "B"+(-8,0)}
{\ar^{} "B"+(8,0); "C"+(-8,0)}
{\ar_{} "D"+(7,0); "E"+(-9,0)}
{\ar^{} "E"+(5,0); "F"+(-11,0)}
{\ar_{a^\ast} "D"+(-3,3); "A"+(-3,-3)}
{\ar^{a^\ast} "E"+(-4,3); "B"+(-4,-3)}
{\ar^{a^\ast} "F"+(-5,3); "C"+(-5,-3)}
{\ar_{} "A"+(-4,-3); "E"+(-9,3)}
{\ar_{} "B"+(-5,-3); "F"+(-10,3)}
\endxy
$

\bigskip

\noindent in which the diagonal morphisms, determined by the Lemma \ref{lemma}, make the upper triangles commute. By Proposition \ref{prop} (4), the colimit of vertical arrows in this diagram is nothing but $a^\ast : \mathbb{L}(B, P) \rightarrow \mathbb{L}(A, P)$. While the colimit of the diagonal morphisms is a right inverse to it. So, $a^\ast$ is a split epi hence a regular epi. 

\end{proof}

Any regular mono $F \rightarrow P$ (or sometimes the object $P$ itself) into a regular $P$ will be called a cover of $F$. Henceforce we assume that a cover is chosen for each left exact presheaf. To prove our main theorems we will use the following two technical lemmas.

\begin{lemma} \label{corollary1}
Every left exact functor $F: \mathbb{C} \rightarrow \mathbb{V}$ is an equalizer of a pair of morphisms between regular functors:

\bigskip

$\xy
(0,0)*{}="A"; 
(16,0)*{}="B";
(34,0)*{}="C";
"A"*{F};
"B"*{P};
"C"*{Q}; 
{\ar^{p} "A"+(2,0); "B"+(-2,0)}
{\ar^{u} "B"+(3,1); "C"+(-3,1)}
{\ar_{v} "B"+(2,-1); "C"+(-4,-1)}
\endxy
$. 

\end{lemma}

\begin{proof}
Let $p: F \rightarrow P$ be the cover for $F$. Being a regular mono $p$ is an equalizer of its cokernel pair. Take $Q$ to be the cover for the cokernel pair of $p$. It is not difficult to see that there is a pair of parallel morphisms between $P$ and $Q$ of which $p$ is an equalizer, as depicted in the diagram above.

\end{proof}

\begin{lemma} \label{corollary2}
Suppose $F$ is a left exact presheaf with the cover $p : F \rightarrow P$, and $T$ is a regular presheaf. There exist a regular presheaf $S$, a regular mono $l: T \rightarrow S$ and a map $w: \mathbb{L}(F, T) \rightarrow \mathbb{L}(P, S)$ such that the following triangle is commutative.

\bigskip

$\xy
(0,0)*{}="A"; 
(36,0)*{}="B";
(18,-18)*{}="C";
"A"*{\mathbb{L}(F, T)};
"B"*{\mathbb{L}(F, S)};  
"C"*{\mathbb{L}(P,S)};    
{\ar^{l_\ast} "A"+(8,0); "B"+(-9,0)}
{\ar_{w} "A"+(3,-3); "C"+(-3,3)}
{\ar_{p^\ast} "C"+(3,3); "B"+(-3,-3)}
\endxy
$

\end{lemma}

\begin{proof}
Let $S$ be determined by the pushout

\bigskip 

$\xy
(0,0)*{}="A"; 
(30,0)*{}="B";
(0,-18)*{}="C";
(30, -18)*{}="D";
"A"*{T}; 
"B"*{S}; 
"C"*{\mathbb{L}(F, T)\cdot F};
"D"*{\mathbb{L}(F, T)\cdot P};
(15,-9)*{p.o.};   
{\ar^{l} "A"+(2,0); "B"+(-2,0)}
{\ar_{\mathbb{L}(F, T)\cdot p} "C"+(10,0); "D"+(-10,0)}
{\ar^{ev} "C"+(-2,3); "A"+(-2,-3)}
{\ar_{\bar{w}} "D"+(-2,3); "B"+(-2,-3)}
\endxy
$

\bigskip
\bigskip

\noindent By transposing $P$ from $\bar{w}$ we get $w: \mathbb{L}(F, T) \rightarrow \mathbb{L}(P, S)$. As shown above in the diagram, the morphism $l$ is a regular mono since it is a pushout of $\mathbb{L}(F, T)\cdot p$, which is a regular mono by Proposition \ref{prop3}. These determine the required data. The commutativity of the triangle is straightforward. 

\end{proof}

Let $\mathbb{R}$ denote the category of regular $\mathbb{V}$-valued functors on $\mathbb{C}$. Of course $\mathbb{R}$ is a subcategory of $\mathbb{L}$.

\begin{theorem}
$\mathbb{R}$ is codense in $\mathbb{L}$.
\end{theorem}

\begin{proof}
$\mathbb{R}$ is codense in $\mathbb{L}$ iff the functor $J: \mathbb{L}^{op} \rightarrow [\mathbb{R}, \mathbb{V}]$ defined by $J(F) = \mathbb{L}(F, -)$ is fully faithful; this means that for each left exact $F$ and $G$, $J_{FG}: \mathbb{L}(G, F) \rightarrow [\mathbb{R}, \mathbb{V}](\mathbb{L}(F, -), \mathbb{L}(G, -))$ is an isomorphism. 

Let us prove this for fixed $F$ and $G$. 

Let $N$ denote $[\mathbb{R}, \mathbb{V}](\mathbb{L}(F, -), \mathbb{L}(G, -)) = \int_X[\mathbb{L}(F, X), \mathbb{L}(G, X)]$.  For a morphism $f: F \rightarrow T$ from $F$ to a regular presheaf, let $\tilde{f}$ be the composite:

\bigskip

$N \cong I\otimes N \stackrel{f\otimes pr_T}{\longrightarrow} \mathbb{L}(F, T)\otimes [\mathbb{L}(F, T), \mathbb{L}(G, T)] \stackrel{ev}{\longrightarrow} \mathbb{L}(G, T)$. 

\bigskip  

\noindent For any map $u: T \rightarrow S$ in $\mathbb{R}$ between regular presheaves we have $\widetilde{uf} = u_\ast \tilde{f}$.

Let $p$, $u$ and $v$ be as in the Lemma \ref{corollary1}. We have: $u_\ast \tilde{p} = \widetilde{up} = \widetilde{vp} = v_\ast \tilde{p}$. Since $p$ is an equalizer of $u$ and $v$ there exists a unique morphism $m: \int_X[\mathbb{L}(F, X), \mathbb{L}(G, X)] \rightarrow \mathbb{L}(G, F)$ such that

\bigskip

$\xy
(0,0)*{}="A"; 
(36,0)*{}="B";
(18,-18)*{}="C";
"A"*{N};
"B"*{\mathbb{L}(G, P)};  
"C"*{\mathbb{L}(G, F)};    
{\ar^{\tilde{p}} "A"+(4,0); "B"+(-8,0)}
{\ar_{m} "A"+(3,-3); "C"+(-3,3)}
{\ar_{p_\ast} "C"+(3,3); "B"+(-3,-3)}
\endxy
$

\bigskip

\noindent is commutative. Let us see that $m$ is a right inverse to $J_{FG}$. All we need is to show that for each $T$ 

\bigskip

$\xy
(0,0)*{}="A"; 
(36,0)*{}="B";
(18,-18)*{}="C";
"A"*{N};
"B"+(6,0)*{[\mathbb{L}(F, T), \mathbb{L}(G, T)]};  
"C"*{\mathbb{L}(G, F)};    
{\ar^{pr_T} "A"+(4,0); "B"+(-10,0)}
{\ar_{m} "A"+(3,-3); "C"+(-3,3)}
{\ar_{\mathbb{L}(-, T)} "C"+(3,3); "B"+(-3,-3)}
\endxy
$

\bigskip

\noindent is commutative. Let $S$, $w$ and $l$ be as in the Lemma \ref{corollary2}. Then we have: 

\bigskip

$N \stackrel{pr_T}{\longrightarrow} [\mathbb{L}(F, T), \mathbb{L}(G, T)] \stackrel{[1, l_\ast]}{\longrightarrow} [\mathbb{L}(F, T), \mathbb{L}(G, S)]$

\bigskip

\noindent equals

\bigskip

$N \stackrel{pr_S}{\longrightarrow} [\mathbb{L}(F, S), \mathbb{L}(G, S)] \stackrel{[l_\ast, 1]}{\longrightarrow} [\mathbb{L}(F, T), \mathbb{L}(G, S)]$ 
 
\bigskip

\noindent equals (using $[l_\ast, 1] = [w, 1][p^\ast, 1]$)  

\bigskip

$N \stackrel{pr_S}{\longrightarrow} [\mathbb{L}(F, S), \mathbb{L}(G, S)] \stackrel{[p^\ast, 1]}{\longrightarrow} [\mathbb{L}(P, S), \mathbb{L}(G,S)] \stackrel{[w, 1]}{\longrightarrow} [\mathbb{L}(F, T), \mathbb{L}(G, S)]$
 
\bigskip

\noindent equals (using $[p^\ast, 1]pr_S = \mathbb{L}(-, S)\tilde{p}$)

\bigskip

$N \stackrel{\tilde{p}}{\longrightarrow} \mathbb{L}(G, P) \stackrel{\mathbb{L}(-, S)}{\longrightarrow} [\mathbb{L}(P, S), \mathbb{L}(G,S)] \stackrel{[w, 1]}{\longrightarrow} [\mathbb{L}(F, T), \mathbb{L}(G, S)]$

\bigskip

\noindent equals (using $\tilde{p} = p_\ast m$)

\bigskip
 
$N \stackrel{m}{\longrightarrow} \mathbb{L}(G, F) \stackrel{p_\ast}{\longrightarrow} \mathbb{L}(G, P) \stackrel{\mathbb{L}(-, S)}{\longrightarrow} [\mathbb{L}(P,S), \mathbb{L}(G, S)] \stackrel{[w, 1]}{\longrightarrow} [\mathbb{L}(F,T), \mathbb{L}(G, S)]$ 
 
\bigskip

\noindent equals (using $[w, 1]\mathbb{L}(-, S)p_\ast = [wp^\ast, 1]\mathbb{L}(-, S) m$)

\bigskip

$N \stackrel{m}{\longrightarrow} \mathbb{L}(G, F) \stackrel{\mathbb{L}(-, S)}{\longrightarrow} [\mathbb{L}(F, S), \mathbb{L}(G, S)] \stackrel{[wp^\ast, 1]}{\longrightarrow} [\mathbb{L}(F,T), \mathbb{L}(G, S)]$

\bigskip

\noindent equals (using $[wp^\ast, 1]\mathbb{L}(-, S) = [l_\ast, 1]\mathbb{L}(-, S) = [1, l_\ast]\mathbb{L}(-, T)$)

\bigskip

$N \stackrel{m}{\longrightarrow} \mathbb{L}(G, F) \stackrel{\mathbb{L}(-, T)}{\longrightarrow} [\mathbb{L}(F, T), \mathbb{L}(G, T)] \stackrel{[1, l_\ast]}{\longrightarrow} [\mathbb{L}(F,T), \mathbb{L}(G, S)]$

\bigskip

\noindent These prove that $[1, l_\ast]pr_T = [1, l_\ast]\mathbb{L}(-, T)m$. Hence $pr_T = \mathbb{L}(-, T)m$ since $[1, l_\ast]$ is a mono.

We have shown that $J_{FG}$ is a split epi. Since $\tilde{p}J_{FG} = p_\ast$ and $p_\ast$ is a mono $J_{FG}$ must be a mono too. To conclude, $J_{FG}$ is an isomorphism. 

\end{proof}

The composition of the Yoneda embedding $Y: \mathbb{C} \rightarrow \mathbb{L}^{op}$ with $J: \mathbb{L}^{op} \rightarrow [\mathbb{R}, \mathbb{V}]$ is a regular fully faithful functor $\mathbb{C} \rightarrow [\mathbb{R}, \mathbb{V}]$. In fact given any subcategory $\mathbb{T}$ of $\mathbb{L}$ all objects of which are regular functors, the canonical evaluation functor  $E: \mathbb{C} \rightarrow [\mathbb{T}, \mathbb{V}]$ is regular. This is because with limits and colimits pointwise in $[\mathbb{T}, \mathbb{V}]$ the evaluation preserves everything that each functor in $\mathbb{T}$ preserves. So in particular, $E$ preserves finite limits and coequalizers of kernel pairs if each presheaf in $\mathbb{T}$ does so. In addition, all we need of $\mathbb{T}$ for $E$ to be fully faithful is a cover of each representable presheaf to be in $\mathbb{T}$ and for each representable $F$ the object $S$ constructed in the Lemma \ref{corollary2} to be in $\mathbb{T}$. 

If $\mathbb{C}$ is a small category (i.e. its set of objects is small), then using simple set theoretic machinery we can find a small subcategory $\mathbb{T}$ of $\mathbb{R}$ with the above properties. Consequently we have:
 
\begin{theorem}
For a small regular $\mathbb{V}$-category $\mathbb{C}$ there exists a small category $\mathbb{T}$ and a regular fully faithful functor $E: \mathbb{C} \rightarrow [\mathbb{T}, \mathbb{V}]$.

\end{theorem}

\end{document}